\magnification=1200

\font\tenbfit=cmbxti10 
\font\sevenbfit=cmbxti10 at 7pt 
\font\sixbfit=cmbxti5 at 6pt 

\newfam\mathboldit 

\textfont\mathboldit=\tenbfit
  \scriptfont\mathboldit=\sevenbfit
   \scriptscriptfont\mathboldit=\sixbfit

\def\bfit           
{\tenbfit           
   \fam\mathboldit 
}

\def\Q{{\bf {Q}}}  \def\F{{\bf F}}
 \def\k{{\bf k}}

\def\Z{{\bf Z}}

\def\Bad{{\bfit Bad}}
\def\oU{{\overline U}}
\def\oV{{\overline V}}
\def\oW{{\overline W}}

\def\house#1{\setbox1=\hbox{$\,#1\,$}%
\dimen1=\ht1 \advance\dimen1 by 2pt \dimen2=\dp1 \advance\dimen2 by 2pt
\setbox1=\hbox{\vrule height\dimen1 depth\dimen2\box1\vrule}%
\setbox1=\vbox{\hrule\box1}%
\advance\dimen1 by .4pt \ht1=\dimen1
\advance\dimen2 by .4pt \dp1=\dimen2 \box1\relax}

  \def\eps{{\varepsilon}}

\def\disp{\displaystyle} 
\def\sm{\smallskip}  \def\noi{\noindent}

\def\build#1_#2^#3{\mathrel{\mathop{\kern 0pt#1}\limits_{#2}^{#3}}}

\def\date {le\ {\the\day}\ \ifcase\month\or janvier
\or fevrier\or mars\or avril\or mai\or juin\or juillet\or
ao\^ut\or septembre\or octobre\or novembre
\or d\'ecembre\fi\ {\oldstyle\the\year}}

\font\fivegoth=eufm5 \font\sevengoth=eufm7 \font\tengoth=eufm10

\newfam\gothfam \scriptscriptfont\gothfam=\fivegoth
\textfont\gothfam=\tengoth \scriptfont\gothfam=\sevengoth

\def\smallsquare{\vbox{\hrule\hbox{\vrule height 1 ex\kern 1 ex\vrule}\hrule}}
\def\cqfd{\hfill \smallsquare\vskip 3mm}


\centerline{}

\vskip 4mm

\centerline{
\bf On the Littlewood conjecture in fields of power series}

\vskip 8mm
\centerline{Boris A{\sevenrm DAMCZEWSKI}  \footnote{}{\rm 
2000 {\it Mathematics Subject Classification : } 11J13, 11J61.} \
\& \ Yann B{\sevenrm UGEAUD} 
\footnote{*}{Supported by the Austrian Science
Fundation FWF, grant M822-N12. } 
}

{\narrower\narrower
\vskip 12mm

\proclaim Abstract. {Let $\k$ be an arbitrary field. 
For any fixed badly approximable power series $\Theta$ in $\k((X^{-1}))$, 
we give an explicit construction of 
continuum many badly approximable 
power series $\Phi$ for which the pair $(\Theta, \Phi)$ satisfies
the Littlewood conjecture. We further discuss the Littlewood
conjecture for pairs of algebraic power series.
}

}

\vskip 6mm
\vskip 8mm

\centerline{\bf 1. Introduction}

\vskip 6mm

A famous problem in simultaneous 
Diophantine approximation is 
the Littlewood conjecture [17]. It claims that, 
for any given pair $(\alpha, \beta)$ of real numbers, we have
$$
\inf_{q \ge 1} \, q \cdot \Vert q \alpha \Vert \cdot \Vert q \beta \Vert = 0,
\eqno (1.1)
$$
where $\Vert \cdot \Vert$ denotes the distance to the nearest integer.
Despite some recent remarkable progress
[24,12], this remains an open problem.

The present Note is devoted to 
the analogous question in fields of power series.
Given an arbitrary field $\k$ and an indeterminate $X$, 
we define a norm $| \cdot |$
on the field $\k ((X^{-1}))$ by setting $|0| = 0$ and,
for any non-zero power series 
$F = F(X) = \sum_{h=-m}^{+ \infty} \, f_h X^{-h}$ with $f_{-m} \not= 0$, 
by setting $|F| = 2^m$. We write $|| F ||$ to denote the norm
of the fractional part of $F$, that is, of the part of the series
which comprises only the negative powers of $X$.
In analogy with (1.1), we ask whether we have
$$
\inf_{q \in \k[X] \setminus \{0\}}  
\, |q| \cdot \Vert q \Theta \Vert \cdot \Vert q \Phi \Vert = 0
\eqno (1.2)
$$
for any given $\Theta$ and $\Phi$ in $\k ((X^{-1}))$.
A negative answer to this question has been obtained by Davenport and
Lewis [11] (see also [3,6,9,10,13] for explicit counter-examples) 
when the field
$\k$ is infinite. As far as we are aware, the problem is still unsolved 
when $\k$ is a finite field (the papers by
Armitage [2], dealing with finite fields of characteristic
greater than or equal to $5$, are erroneous, as kindly pointed out to us
by Bernard de Mathan). 

A first natural question regarding this problem can be stated as follows:
\medskip

\noindent{\it Question 1.} Given 
a badly approximable power series $\Theta$, 
does there exist a power series
$\Phi$ such that the pair $(\Theta, \Phi)$ satisfies 
{\it non-trivially} the Littlewood conjecture?
\medskip

\noindent First, we need to explain what is meant by {\it non-trivially} 
and why we restrict our attention to badly approximable power series, that
is, to power series from the set 
$$
\Bad = \bigl\{ \Theta \in \k ((X^{-1})) :
\inf_{q \in \k[X] \setminus \{0\}}  \,
|q| \cdot \Vert q \Theta \Vert > 0 \bigr\}.
$$
Obviously, (1.2) holds as soon as $\Theta$ or $\Phi$ does not belong to $\Bad$.
This is also the case when $1$, $\Theta$ and $\Phi$ are linearly dependent 
over $\k[X]$. Hence, by {\it non-trivially}, we simply mean that both of these
cases are excluded.

\medskip

In the present paper, we answer positively Question 1 by
using the constructive approach developed in [1]. 
Our method rests on the basic theory of continued fractions and 
works without any restriction on the field $\k$. 
Actually, our result is much more precise and  motivates the investigation 
of a stronger question, introduced and discussed in Section 2. 
Section 3 is concerned with the Littlewood conjecture for pairs
of algebraic power series. When $\k$ is a finite field,
we provide several examples of such pairs for which (1.2) holds.  
In particular, we show that there exist infinitely many pairs of 
quartic power series in $\F_3((X^{-1}))$ that
satisfy non-trivially the Littlewood conjecture. 
It seems that no such pair was previously known. 
The proofs are postponed to Sections 5 and 6, 
after some preliminaries on continued fractions gathered in Section 4.

\vskip 8mm

\centerline{\bf 2. Main results}

\vskip 6mm

The real analogue of Question 1 was answered positively 
by Pollington and Velani [24] by
using metric theory of Diophantine approximation, as a consequence
of a much stronger statement. Some years later, 
Einsiedler, Katok and Lindenstrauss [12] proved the
outstanding result that the set of pairs of real numbers for which
the Littlewood conjecture does not hold has Hausdorff dimension zero.
Obviously, this implies a positive answer to Question 1. 
However, it is unclear that either of these methods could be transposed in 
the power series case. Furthermore, both methods are not constructive,
in the sense that they do not yield
explicit examples of pairs of real numbers satisfying (1.1).

A new, explicit and elementary approach to solve the real
analogue of Question 1 is developed in [1].
It heavily rests on the theory of continued fractions
and it can be quite naturally adapted to the
function field case. Actually,  
our Theorem 1 answers a strong form of Question 1.

\proclaim Theorem 1.
Let $\varphi$ be a positive, non-increasing function defined on the set of
positive integers and with $\varphi(1) = 1$ and
$\lim_{d \to + \infty} \, \varphi(d) = 0$.
Given $\Theta$ in $\Bad$, there is an uncountable subset
$B_{\varphi}(\Theta)$ of $\Bad$ such that,
for any $\Phi$ in $B_{\varphi}(\Theta)$, the power series
$1$, $\Theta$, $\Phi$ are linearly independent over
$\k[X]$ and there exist 
polynomials $q$ in $\k[X]$ with arbitrarily large degree and satisfying
$$
|q|^2 \cdot \Vert q \Theta \Vert \cdot \Vert q \Phi \Vert \le 
{1 \over \varphi(|q|)}. \eqno (2.1)
$$
In particular, the Littlewood conjecture holds non-trivially for 
the pair $(\Theta, \Phi)$ for any $\Phi$ in $B_{\varphi}(\Theta)$. 
Furthermore, the set $B_{\varphi}(\Theta)$ can be effectively
constructed. 

Although the proof closely follows that of Theorem 1
from [1], we give it with full detail. 
Actually, some steps are even slightly easier than in the real case. 

\medskip

Observe that, for any
given $\Theta$ and $\Phi$ in $\Bad$, 
there exists a positive constant
$c(\Theta, \Phi)$ such that
$$
|q|^2 \cdot \Vert q \Theta \Vert \cdot \Vert q \Phi \Vert \ge
c(\Theta, \Phi)
$$
holds for any non-zero polynomial $q$ in $\k[X]$. 
In view of this
and of Theorem 1, we propose the following problem in which we ask 
whether the above inequality is best possible.

\medskip

\noindent{\it Question 2.} Given 
a power series $\Theta$ in $\Bad$, does there exist a power series
$\Phi$ such that the pair $(\Theta, \Phi)$ satisfies 
{\it non-trivially} the Littlewood conjecture and such that we
moreover have
$$
\liminf_{\deg q \to + \infty} \,
|q|^2 \cdot \Vert q \Theta \Vert \cdot \Vert q \Phi \Vert < + \infty \; ?
\eqno (2.2)
$$
\medskip

\noindent  The restriction `non-trivially' in the statement
of Question 2 is needed, since (2.2) clearly holds when
the power series $1$, $\Theta$, $\Phi$ are linearly dependent over
$\k[X]$. There are, however, non-trivial examples
for which (2.2) holds.
Indeed, if the continued fraction expansion of a power series $\Theta$ 
begins with 
infinitely many palindromes and if $\Phi = 1/\Theta$, then (2.2) 
is true for the pair $(\Theta, \Phi)$.
This can be seen by working out in the power series case the
arguments from Section 4 of [1].

\proclaim Theorem 2. 
Let $\Theta$ be an element of the field $\k ((X^{-1}))$ whose 
continued fraction expansion
begins with infinitely many palindromes. 
Then, the Littlewood conjecture is true
for the pair $(\Theta, \Theta^{-1})$ and, furthermore, we have
$$
\liminf_{\deg q \to + \infty} \,
|q|^2 \cdot \Vert q \Theta \Vert \cdot \Vert q \Theta^{-1} \Vert < + \infty.
$$
Moreover, if $\k$ has characteristic zero, then $\Theta$ is
transcendental over $\k (X)$.

We can weaken the assumption that the continued fraction 
expansion of $\Theta$ begins with infinitely many palindromes 
to get additional examples of pairs $(\Theta,\Theta^{-1})$ 
that satisfy the Littlewood conjecture.
Before stating our next result, we need
to introduce some notation. It is convenient to use the terminology
from combinatorics on words.
We identify any finite or infinite word
$W = w_1 w_2 \ldots$ on the alphabet $\k[X] \setminus \k$ 
with the sequence of partial quotients $w_1, w_2, \ldots $
Further, if $U = u_1 \ldots u_m$ and 
$V = v_1 v_2 \ldots$ are words on $\k[X] \setminus \k$, 
with $V$ finite or infinite, and if $u_0$ is in $\k[X]$,
then $[u_0, U, V]$ denotes the continued fraction 
$[u_0, u_1, \ldots, u_m, v_1, v_2,\ldots]$. 
The mirror image of any finite word $W = w_1 \ldots w_m$
is denoted by $\oW := w_m \ldots w_1$.
Recall that a
palindrome is a finite word $W$ such that $\oW=W$.
Furthermore, we denote by $|W|$ the number
of letters composing $W$ (here, we clearly have $|W| = m$). There
should not be any confusion between $|W|$ and the norm $|F|$ of
the power series $F$.

\proclaim Theorem 3. 
Let $\Theta$ be in $\Bad$. Denote 
by $(p_n/q_n)_{n \ge 1}$ the sequence of its convergents.
Assume that there exist a positive real number $x$, 
a sequence of finite words $(U_k)_{k \ge 1}$, 
and a sequence of palindromes $(V_k)_{k \ge 1}$
such that, for every $k \ge 1$, the continued fraction expansion of
$\Theta$ is equal to $[U_k, V_k\ldots]$ 
and $|V_{k+1}| > |V_k| \ge x |U_k|$.
Set further 
$$
M = \limsup_{\ell \to + \infty} \, {\deg q_{\ell} \over \ell}
\quad {\it and} \quad
m = \liminf_{\ell \to + \infty} \, {\deg q_{\ell} \over \ell}.
$$
If we have
$$
x > 3\cdot{M \over m}-1, \eqno (2.3)
$$
then the Littlewood conjecture is true for the pair $(\Theta, \Theta^{-1})$.
Moreover, if $\k$ has characteristic zero, then $\Theta$ is
transcendental over $\k (X)$.

From now on, we make use of the following notation: 
if $\ell$ is a positive integer, then $W^{[\ell]}$
denotes the word obtained by concatenation of $\ell$ copies of the word $W$. 

\proclaim Theorem 4. 
Let $\Theta = [a_0, a_1, a_2, \ldots]$ be in $\Bad$. Denote 
by $(p_n/q_n)_{n \ge 1}$ the sequence of its convergents.
Assume that there exist a finite word $V$, 
a sequence of finite words $(U_k)_{k \ge 1}$, 
an increasing sequence of positive integers  $(n_k)_{k \ge 1}$ 
and a positive real number $x$ 
such that, for every $k \ge 1$, the continued fraction expansion of
$\Theta$ is equal to $[U_k, V^{[n_k]}\ldots]$ and $|V^{[n_k]}| \ge  x|U_k|$. 
Let $\Phi$ be the quadratic power series defined by  
$$
\Phi:=[\oV,\oV,\oV,\ldots].
$$
Set further 
$$
M = \limsup_{\ell\to+\infty} \deg a_{\ell}
\quad {\it and} \quad
m = \liminf_{\ell \to + \infty} \deg a_{\ell}.
$$
If we have
$$
x > {M \over m}, \eqno (2.4)
$$
then the pair $(\Theta, \Phi)$ satisfies the Littlewood conjecture.
Moreover, if $\k$ has characteristic zero, then $\Theta$ is
transcendental over $\k (X)$.

The last assertion of Theorems 2, 3 and 4 
follows from the analogue of
the Schmidt Subspace Theorem in fields of power series over a field
of characteristic zero, worked out by Ratliff [25]. 
It is well-known that the analogue of the Roth theorem 
(and, {\it a fortiori}, the analogue of the Schmidt Subspace Theorem)
does not hold for fields of power series over a finite field. 
For $\k = \F_p$, a celebrated example given by Mahler [18]
is recalled in Section 3. 

\medskip

Theorems 2, 3 and 4 will be used in the next section 
to provide new examples of pairs of algebraic power series 
satisfying the Littlewood conjecture.
\vskip 8mm
\goodbreak

\centerline{\bf 3. On the Littlewood conjecture for pairs
of algebraic power series}

\vskip 6mm

It is of particular interest to determine whether
the Littlewood conjecture holds for pairs of algebraic real numbers. 
To the best of our knowledge, only two results are known 
in this direction. 
First, if $(\alpha,\beta)$ is a pair of real numbers lying 
in a same quadratic field, then $1$, $\alpha$ and $\beta$ are 
linearly dependent over $\Q$ and the Littlewood conjecture is 
thus easily satisfied. This was for instance remarked in [7].  
The other result is due to Cassels and Swinnerton-Dyer [8] 
who proved that the Littlewood conjecture is satisfied for pairs 
of real numbers lying in a same cubic field. However, it is generally believed 
that no algebraic number of degree greater than or
equal to $3$ is badly approximable. 
At present, no pair of algebraic numbers 
is known to satisfy non-trivially the Littlewood conjecture.

\medskip

In this Section, we discuss whether the 
(function field analogue of the) Littlewood conjecture 
holds for pairs of algebraic power series defined over a finite field $\k$. 
Our knowledge is slightly better than in the real case, 
especially thanks to works of Baum and Sweet [4] 
and of de Mathan [19,20,21,22] that we recall below. 

First, we observe that, as in the real case, (1.2) holds when 
$\Theta$ and $\Phi$ are in a same
quadratic extension of $\k [X]$, since $1$, $\Theta$ and $\Phi$ are then
linearly dependent over $\k [X]$. We further observe that the 
existence of the Frobenius automorphism (that is, the
$p$-th power map) yield   
many examples of well-approximated algebraic power series. 
For instance, for any prime number $p$,
the power series $\Theta_p = [0,X,X^p,X^{p^2},X^{p^3},\ldots]$ 
is a root in $\F_p((X^{-1}))$ of the polynomial $Z^{p+1} + X Z -1$,
and $\Theta_p$ is well-approximated 
by rational functions. Indeed, there exist infinitely 
many rational functions $p_n/q_n$ such that 
$$
\biggl\vert \Theta_p-{p_n\over q_n}\biggr\vert\le 
{1\over \vert q_n\vert^{p+1}} \cdot
$$
Clearly, for any (algebraic or transcendental)
power series $\Phi$ in $\F_p((X^{-1}))$,
the Littlewood conjecture holds for the pair $(\Theta_p, \Phi)$.

On the other hand, there are several results on pairs 
of algebraic functions that satisfy non-trivially the Littlewood conjecture.
De Mathan [21]
established that (1.2) holds for any pair $(\Theta, \Phi)$ of
quadratic elements when $\k$ is any finite 
field of characteristic $2$ (see also 
[19,20] for results when $\k$ is any finite field). 
Furthermore, he proved in [22] the analogue of 
the Cassels and Swinnerton-Dyer theorem when $\k$ is a finite field. 
We stress that, when $\k$ is finite, there do exist,
unlike in the real case, algebraic power series in $\Bad$ that are of
degree greater than or equal to $3$ over $\k(X)$.
The first example was given by Baum and Sweet [4] who proved that,
for $\k = \F_2$, the unique $\Theta$ in $\F_2 ((X^{-1}))$ which satisfies
$X \Theta^3 + \Theta + X = 0$ is in $\Bad$. 
Thus, it follows from [22] that the pair of algebraic 
power series $(\Theta,\Theta^{-1})$ satifies non-trivially 
the Littlewood conjecture.  

Further examples of badly approximable algebraic power series 
were subsequently found by several authors. 
It turns out that many of these examples contain some   
symmetric patterns in their continued fraction expansion. 
In the sequel of this Section, this property is used in
order to apply Theorems 2, 3 and 4 
to provide new examples of pairs of algebraic power series 
satisfying non-trivially the Littlewood conjecture.
These examples also illustrate the well-known fact that there
is no analogue to the
Schmidt Subspace Theorem for power series over finite fields.

We keep on using the terminology from combinatorics on words.
For sake of readability we sometimes write commas to separate the 
letters of the words we consider.

\medskip

\vskip 6mm

{\bf 3.1. A first example of a badly approximable 
quartic in  $\F_3((X^{-1}))$}

\vskip 6mm

Mills and Robbins [23] established that the polynomial 
$$
X(X+2)Z^4 - (X^3 + 2 X^2 + 2 X + 2) Z^3 + Z - X - 1
$$
has a root $\Theta$ in $\F_3((X^{-1}))$ 
whose continued fraction expansion is expressed as follows.
For every positive integer $n$, set 
$$
H_n=X^{[3^n-2]},X+\eps,2X+\eps, (2X)^{[3^n-2]},2X+\eps,X+\eps,
$$
where $\eps=2$ if $n$ is odd and $\eps=1$ otherwise. 
Then, the continued fraction expansion of 
the quartic power series $\Theta$ is given by 
$$
\Theta=[X,2X+2,X+1,H_1,H_2,H_3,\ldots].
$$

\medskip

It turns out that the continued fraction expansion of $\Theta$ 
contains some symmetric patterns that we can use to apply Theorem 3. 
This gives rise to the following result.

\proclaim Theorem 5. The pair $(\Theta, \Theta^{-1})$ satisfies 
the Littlewood conjecture. In particular, there exists a pair of quartic
power series in $\F_3((X^{-1}))$ satisfying non-trivially 
the Littlewood conjecture. 

To our knowledge, this is the first known example of a pair of algebraic
power series of degree greater than $3$ for which the Littlewood
conjecture is non-trivially satisfied.

\medskip

\noi {\bf Proof.} 
For every integer $n\ge 2$, set
$$
U_n := X,2X+2,X+1,H_1,H_2,H_3, \ldots , H_{n-1}
$$
and
$$
V_n := H_n, X^{[3^n-2]}.
$$
Since $X^{[3^n-2]}$ is a prefix of $H_{n+1}$, 
the continued fraction expansion of $\Theta$ is equal
to $[U_nV_n\ldots]$. Furthermore, $V_n$ is a palindrome and 
the length of $H_n$ (resp. of $U_n$, of $V_n$)
is equal to $2\cdot3^n$ (resp. to $3^n$, to $3^{n+1}-2$). 
In particular,  we have $|V_n| >2.5|U_n|$ for every $n\ge 2$, and, 
since all the partial quotients of $\Theta$ are linear, 
the assumption (2.3) is satisfied.
We apply Theorem 3 to complete the proof. 
\cqfd
\medskip

\vskip 6mm
\goodbreak

{\bf 3.2. An infinite family of badly approximable   
quartics in $\F_3((X^{-1}))$}

\vskip 6mm

We now consider the infinite family of badly 
approximable quartics in $\F_3((X^{-1}))$ 
introduced by Lasjaunias in [15]. 
Let $k$ be a non-negative integer. For any non-negative integer $n$, set
$$
u_n = (k+2) 3^n - 2,
$$
and define the finite word $H_n(X)$ on $\F_3 [X] \setminus \F_3$ by
$$
H_n(X) := (X+1) X^{[u_n]} (X+1).
$$
Then, consider the power series
$$
\Theta(k) := [0, H_0 (X), H_1 (-X), H_2 (X), \ldots , 
H_n ((-1)^n X)), \ldots]. \eqno (3.1)
$$
This definition obviously implies that $\Theta(k)$ 
is badly approximable by rational functions, 
since all of its partial quotients are linear.
Lasjaunias [15] established that $\Theta(k)$ 
is a quartic power series. 
More precisely,  if $(p_n(k)/q_n(k))_{n\ge 0}$ 
denotes the sequence of convergents to $\Theta(k)$, he proved that  
$$
q_k(k)\Theta^4(k)-p_k(k)\Theta^3(k)+q_{k+3}(k)\Theta(k)-p_{k+3}(k)=0.
$$
The description of the continued fraction expansion 
of $\Theta(k)$ given in (3.1) makes transparent the 
occurrences of some palindromic patterns. 
This can be used to apply Theorem 3 and  yields the following result.

\proclaim Theorem 6. For any non-negative integer $k$, 
the pair $(\Theta(k), \Theta(k)^{-1})$ satisfies the Littlewood conjecture. 
In particular, there exist infinitely many pairs of quartic
power series in $\F_3((X^{-1}))$ 
satisfying non-trivially the Littlewood conjecture.

\medskip

\noi {\bf Proof.} 
For any even positive integer $n$, set
$$
U_n := H_0 (X) H_1 (-X) H_2(X)\ldots H_{n-2} (X) (-X+1)
$$
and
$$
V_n := (-X)^{[u_{n-1}]} (-X+1) (X+1) X^{[u_n]} (X+1) (-X+1) 
(-X)^{[u_{n-1}]}.
$$
Observe that the continued fraction expansion of $\Theta (k)$ is equal
to $[0,U_n V_n\ldots]$ and that
$$
|U_n| = \biggl( {k+2 \over 2} \biggr) 3^{n-1} - {k \over 2}   \quad
{\rm and} \quad |V_n| = 5 (k+2) 3^{n-1} - 2.
$$
Furthermore, $V_n$ is a palindrome.
We have $|V_n| \ge 3|U_n|+3$ for $n \ge 2$, and, 
since all the partial quotients of $\Theta (k)$ are linear, 
the assumption (2.3) is satisfied.
We apply Theorem 3 to complete the proof. 
\cqfd

\medskip

\vskip 6mm

{\bf 3.3. Badly approximable power series in $\F_p((X^{-1}))$ with $p\ge 5$}

\vskip 6mm

Let $p \ge 5$ be a prime number. For any positive integer $k$, 
consider the polynomial $f_k$ in $\F_p[X]$ defined by  
$$
f_k=\sum {k-j \choose j}X^{k-2j},
$$
where the sum is over all integers $j$ such that $0\le 2j\le k$. 
Then, Mills and Robbins [23] showed that   
the polynomial of degree $p+1$  
$$
X Z^{p+1} - (X^2 - 3) Z^p + (X f_{p-2} - 3 f_{p-1}) Z
- f_{p-2} (X^2 - 3) + f_{p-1} X
$$
has a root $\Theta_p$ in $\F_p((X^{-1}))$ 
with a nice continued fraction expansion. 
Let $V(-1)=-X,-X$ and $V(3)=X/3,3X$ and, for $k \ge 1$,
set 
$$
L_k(-1)=V(-1)^{[(p^k-1)/2]} \quad
{\rm and} \quad L_k(3)=V(3)^{[(p^k-1)/2]}.
$$ 
Mills and Robbins proved that the continued 
fraction expansion of $\Theta_p$ is given by 
$$
\Theta_p=[X,L_0(3), -X/3,L_0(-1),X,L_1(3), -X/3,L_1(-1),
X,L_2(3), -X/3,L_2(-1),\ldots],
$$
where $L_0(3)$ and $L_0(-1)$ are equal to the empty word.
It follows that $\Theta_p$ is badly approxi\-mable by rational functions, 
all of its partial quotients being linear. Moreover, 
$\Theta_p$ is not quadratic since its continued 
fraction expansion is not eventually periodic. 

As a consequence of Theorem 3, we get the following result.

\proclaim Theorem 7. For any prime number $p \ge 7$, 
the pair $(\Theta_p, \Theta_p^{-1})$ of algebraic power series  
in $\F_p((X^{-1}))$ satisfies non-trivially the Littlewood conjecture.  

Moreover, we can apply Theorem 4 to provide pairs 
of algebraic power series of distinct degrees 
satisfying non-trivially the Littlewood conjecture. 
To the best of our knowledge, no such pair was previously known.

\proclaim Theorem 8. 
Let $p \ge 5$ be a prime number. Let $\Theta_p$ be as above.
Let $\Phi_p$ be the quadratic power series in $\F_p((X^{-1}))$ defined by 
$$ 
\Phi_p:=[3X,X/3,3X,X/3,3X,X/3,3X,\ldots].
$$
Then the pair $(\Theta_p, \Phi_p)$ satisfies non-trivially 
the Littlewood conjecture.  

\medskip

\noi {\bf Proof of Theorems 7 and 8.} 
For any even positive integer $n$, set
$$
U_n := X,-X/3,X,L_1(3), 
-X/3,L_1(-1),X, L_2(3), -X/3,L_2(-1), X \ldots ,L_{n-1}(-1),X
$$
and
$$
V_n := (X/3,3X)^{[(p^n-3)/2]}, X/3.
$$
Observe that the continued fraction expansion of $\Theta_p$ is
equal to $[U_nV_n\ldots]$ and that
$$
|U_n| = 1+\left(2\cdot {p^n-1\over p-1}\right)  \quad
{\rm and} \quad |V_n| = p^n-2.
$$
Furthermore, $V_n$ is a palindrome and $|V_n| \ge 2.5|U_n|$ 
holds for $p \ge 7$ and $n \ge 2$.  
Since all the partial quotients of $\Theta_p$ are linear, 
the assumption (2.3) is then satisfied.
We apply Theorem 3 to complete the proof of Theorem 7. 

To get Theorem 8, we observe that $L_n (3)$ is the concatenation
of $(p^n-1)/2$ copies of the word $V(3)$, and we check that
$|L_n(3)| \ge 1.5|U_n|$ holds for $p\ge 5$ and $n \ge 2$.  
Since all the partial quotients of $\Theta_p$ are linear, 
the assumption (2.4) is then satisfied.
We then apply Theorem 4 to complete the proof of Theorem 8.  \cqfd

\medskip

\vskip 6mm

{\bf 3.4. A normally approximable quartic in $\F_3((X^{-1}))$}

\vskip 6mm

We end this Section with another quartic power series in $\F_3((X^{-1}))$ 
found by Mills and Robbins [23].  
Unlike the previous examples, this quartic is not badly 
approximable but we will see that it has some 
interesting Diophantine properties. 
 
Mills and Robbins pointed out that the polynomial 
$$
Z^4 + Z^2 - X Z + 1 
$$
has a unique root $\Theta$ in $\F_3 ((X^{-1}))$. 
 They observed empirically that 
$\Theta$ has a particularly simple continued fraction expansion.
Define recursively a sequence $(\Omega_n)_{n \ge 0}$ of words on
the alphabet $\F_3[X] \setminus \F_3$ by setting $\Omega_0 = \varepsilon$,
the empty word, $\Omega_1 = X$, and
$$
\Omega_n = \Omega_{n-1} (-X) \Omega_{n-2}^{(3)} (-X) \Omega_{n-1}, \quad
\hbox{for $n \ge 2$}. \eqno (3.2)
$$
\noindent Here, if $W=w_1w_2\ldots w_r=w_1,w_2,\ldots,w_r$ with 
$w_i\in \F_3[X]\setminus \F_3$, then $W^{(3)}$ denotes the word obtained by 
taking the cube of every letter of $W$, that is,
$W^{(3)}:=w_1^3,w_2^3,\ldots,w_r^3$. Set
$$
\Omega_{\infty} = \lim_{n \to + \infty} \, \Omega_n. \eqno (3.3)
$$
Buck and Robbins [5] confirmed a conjecture of Mills and
Robbins [23] asserting that the continued
fraction expansion of $\Theta$ is $[0, \Omega_{\infty}]$ (note that their
proof was later simplified by Schmidt [26], 
and that Lasjaunias [14] gave an alternative proof).  

The quartic power series $\Theta$ does not lie in $\Bad$. 
Lasjaunias [14, Theorem A] proved that $\Theta$ is normally approximable 
(this terminology is explained in [16]) in the following sense:
there exist infinitely many rational functions $p/q$ such that 
$$
\vert \Theta-p/q\vert\le \vert q\vert^{-(2+2/\sqrt{3\deg q})},
$$ 
while for any positive real number $\eps$ there are only finitely 
many rational functions $p/q$ such that 
$$
\vert \Theta-p/q\vert\le \vert q\vert^{-(2+2/\sqrt{3\deg q}+\eps)}.
$$ 

Note that an easy induction based on (3.2) shows that 
for every positive integer $n$ the word $\Omega_n$ is a palindrome.   
By (3.3), we thus get that the continued fraction expansion of $\Theta^{-1}$ 
begins with infinitely many palindromes. 
The following consequence of Theorem 2 and of Theorem A from [14] is
worth to be pointed out.

\proclaim Theorem 9. Let $\Theta$ be the unique root in
$\F_3((X^{-1}))$ of the polynomial 
$Z^4+Z^2-X Z+1$.
Then, 
$$
\inf_{q \in \k[X] \setminus \{0\}}  
\, |q|^2 \cdot \Vert q \Theta \Vert \cdot \Vert q \Theta^{-1} \Vert < + \infty
$$
and for any positive real number $\eps$ we have 
$$
|q|^{2 + 4/ \sqrt{3\deg q}+\eps}
\cdot \Vert q \Theta \Vert \cdot \Vert q \Theta^{-1} \Vert \ge 1,
$$
for any $q$ in $\F_3 [X]$ with $\deg q$ large enough.

\vskip 8mm

\centerline{\bf 4. Preliminaries on continued fraction expansions
of power series}

\vskip 6mm

It is well-known that the continued fraction algorithm can as well be applied
to power series. The partial quotients are then
elements of $\k[X]$ of positive degree.
We content ourselves to recall some basic facts, and we direct the
reader to Schmidt's paper [26] and to Chapter 9 of Thakur's book [27] for
more information. 

Specifically, given a power series
$F = F(X)$ in $\k((X^{-1}))$, which we assume not to be a rational function,
we define inductively the sequences $(F_n)_{n \ge 0}$
and $(a_n)_{n \ge 0}$ by $F_0 = F$ and
$F_{n+1} = 1/(F_n - a_n)$, where $a_n = \Vert F_n \Vert$.
Plainly, for $n \ge 1$,
the $a_n$ are polynomials of degree at least one.
We then have 
$$
F = [a_0 , a_1, a_2, \ldots] a_0 + {1 \over \disp a_1 +
{\strut 1 \over \disp a_2 + {\strut 1 \over \disp \ldots}}} \cdot
$$
The truncations $[a_0 , a_1, a_2, \ldots , a_n]:= p_n / q_n$, with relatively
prime polynomials $p_n$ and $q_n$, are rational functions and are 
called the {\it convergents} to $F$.
It is easily seen that
$$
\deg q_{n+1} = \deg a_{n+1} + \deg q_n,
$$
thus
$$
\deg q_n = \sum_{j = 1}^n \, \deg a_j. \eqno (4.1)
$$
Furthermore, we have
$$
\deg (q_n F - p_n) = - \deg q_{n+1} < - \deg q_n,
$$
that is,
$$
\Vert q_n F \Vert = | q_{n+1} |^{-1}
= 2^{-\deg q_{n+1}} < 2^{-\deg q_n}. \eqno(4.2)
$$

We stress that $F$ is in $\Bad$ if and only if the degrees of
the polynomials $a_n$ are uniformly bounded. 
We also point out that $| \cdot |$ is an
ultrametric norm, that is, $|F + G| \le \max\{|F|, |G|\}$
holds for any $F$ and $G$ in $\k ((X^{-1}))$, with equality 
if $|F| \not= |G|$. 

We end this Section by stating three basic lemmas on 
continued fractions in $\k((X^{-1}))$.

\proclaim Lemma 1.
Let $\Theta = [a_0, a_1, a_2, \ldots]$ be an element of $\k((X^{-1}))$  
and let
$(p_n/q_n)_{n \ge 1}$ be its convergents. Then, for any $n \ge 2$, we have
$$
{q_{n-1} \over q_n } = [0, a_n, a_{n-1}, \ldots , a_1].
$$

\medskip

\noindent {\bf Proof.} As in the real case, 
this easily follows from the recursion formula 
$q_{n+1} = a_{n+1} q_n + q_{n-1}$. \cqfd

\medskip

\proclaim Lemma 2.
Let $\Theta = [0, a_1, a_2, \ldots]$ and $\Phi [0, b_1, b_2, \ldots]$ be two elements of $\k((X^{-1}))$. 
Assume that there exists a positive integer $n$ such that
$a_i = b_i$ for any $i=1, \ldots, n$. We then have
$\vert\Theta - \Phi\vert \le \vert q_{n} \vert^{-2}$, 
where $q_{n}$ denotes the denominator
of the $n$-th convergent to~$\Theta$.

\medskip

\noindent {\bf Proof.} Let $p_n/q_n$ be the $n$-th convergent to $\Theta$. 
By assumption, $p_n/q_n$ is also the $n$-th convergent 
to $\Phi$ and we have 
$$
\vert\Theta - \Phi\vert \le \max\left\{\vert\Theta-p_n/q_n\vert, \;  
\vert\Phi-p_n/q_n\vert\right\}\le \vert q_{n} \vert^{-2},
$$
since the norm $\vert\cdot\vert$ is ultrametric. \cqfd

\medskip

\proclaim Lemma 3. Let $M$ be a positive real number. Let 
$\Theta = [0, a_1, a_2, \ldots]$ and $\Phi [0, b_1, b_2, \ldots]$ be two elements of $\k((X^{-1}))$ 
whose partial quotients are of degree at most $M$. 
Assume that there exists a positive integer $n$ such that  
$a_i = b_i$ for any $i=1, \ldots, n$ and $a_{n+1}\not=b_{n+1}$. 
Then, we have
$$
|\Theta - \Phi| \ge {1\over 2^{2M} \vert q_n\vert^2},
$$ 
where $q_n$ denotes the denominator
of the $n$-th convergent to $\Theta$.

\medskip

\noindent {\bf Proof.} Set $\Theta'=[a_{n+1}, a_{n+2},\ldots]$ and 
$\Phi'=[b_{n+1}, b_{n+2},\ldots]$. Since $a_{n+1}\not= b_{n+1}$, we have 
$$
\vert\Theta'-\Phi'\vert\ge 1. \eqno (4.3)  
$$
Furthermore, since the degrees of the partial quotients 
of both $\Theta$ and $\Phi$ are
bounded by $M$, we immediately obtain that
$$
\vert\Theta'\vert \le 2^M \quad{\rm{and}} \quad 
\vert\Phi'\vert\le 2^M. \eqno (4.4)
$$
Denote by $(p_j / q_j)_{j \ge 1}$ 
the sequence of convergents to $\Theta$. Then, the theory of continued
fractions gives that 
$$
\Theta={p_n\Theta'+p_{n-1}\over q_n\Theta'+q_{n-1}} \quad {\rm and} \quad
\Phi={p_n\Phi'+p_{n-1}\over q_n\Phi'+q_{n-1}},
$$
since the first $n$-th partial quotients of $\Theta$ and
$\Phi$ are assumed to be the same.  
We thus obtain 
$$
\vert\Theta-\Phi\vert\left\vert{p_n\Theta'+p_{n-1}\over q_n\Theta'+q_{n-1}}-
{p_n\Phi'+p_{n-1}\over q_n\Phi'+q_{n-1}}\right\vert
= \left\vert{ \Theta'-\Phi'
\over (q_n\Theta'+q_{n-1})(q_n\Phi'+q_{n-1})}
\right\vert=\left\vert{ \Theta'-\Phi'
\over \Theta'\Phi'q_n^2}\right\vert\cdot
$$
Together with (4.3) and (4.4), this yields  
$$
\vert\Theta-\Phi\vert\ge {1\over 2^{2M} \vert q_n\vert^2},
$$
concluding the proof of the lemma. 
\cqfd

\vskip 8mm

\centerline{\bf 5. Proof of Theorem 1}

\vskip 6mm

Without any loss of generality, we may assume that $|\Theta| \le 1/2$ and we
write $\Theta =[0, a_1, a_2, \ldots, a_k, \ldots]$.
Let $M$ be an upper bound for the degrees of the polynomials $a_k$.
We first construct inductively a rapidly increasing sequence
$(n_j)_{j \ge 1}$ of positive integers.  
We set $n_1 = 1$ and we proceed with the inductive step.
Assume that $j \ge 2$ is such that $n_1, \ldots , n_{j-1}$ have been
constructed. Then, we choose $n_j$ sufficiently large in order that  
$$
\varphi(2^{m_j}) \le 2^{-2 (M+2)(m_{j-1}+1)}, \eqno (5.1)
$$
where $m_j=n_1+n_2+\ldots +n_j+(j-1)$. Such a choice  
is always possible since $\varphi$ tends to zero at infinity and since
the right-hand side of (5.1) only depends on $n_1,n_2,\ldots,n_{j-1}$.

Our sequence $(n_j)_{j \ge 1}$ being now constructed, 
for an arbitrary sequence  ${\bf t}=(t_j )_{j \ge 1}$ 
with values in $\k[X] \setminus \k$, we set 
$$
\eqalign{
\Phi_{\bf t} = & [0, b_1, b_2, \ldots] \cr
= & [0, a_{n_1}, \ldots, a_1, t_1, 
a_{n_2}, \ldots, a_1, t_2, 
a_{n_3}, \ldots, a_1, t_3,\ldots ]. \cr}
$$
Then, we introduce the set 
$$
B_{\varphi}(\Theta)=\left\{\Phi_{\bf t}, \; {\bf t}\in (\k_{M+1}[X] \cup
\k_{M+2} [X])^{\Z_{\ge 1}}\right\},
$$
where $\k_n [X]$ denotes the set of polynomials in $\k[X]$ of
degree $n$.
Clearly, the set $B_{\varphi}(\Theta)$ is uncountable.

Let $\Phi$ be in $B_{\varphi}(\Theta)$. We first prove
that (2.1) holds for the pair $(\Theta, \Phi)$. 
Denote by $(p_n / q_n )_{n \ge 1}$ 
(resp. by $(r_n  / s_n )_{n \ge 1}$)
the sequence of convergents to $\Theta$ (resp. to $\Phi$). 
Let $j \ge 2$ be an integer.
We infer from Lemma 1 that
$$
{s_{m_j - 1}  \over s_{m_j} } = [0, a_1, \ldots , a_{n_j}, 
t_{j-1}, a_1, \ldots, a_{n_{j-1}}, t_{j-2} ,
\ldots, t_1, a_1, \ldots, a_{n_1} ].
$$
By (4.2), we have   
$$
\Vert s_{m_j}  \, \Phi \Vert
\le \vert s_{m_j}\vert^{-1}. \eqno (5.2)
$$
On the other hand, Lemma 2 implies that 
$$
\biggl\vert \Theta - {s_{m_j - 1}  \over s_{m_j} } 
\biggr\vert
\le {1\over \vert q_{n_j}\vert^2} = 2^{- 2 \deg q_{n_j} }. 
$$
Consequently, we get
$$
\Vert s_{m_j}  \, \Theta \Vert
\le 2^{\deg s_{m_j}  - 2 \deg q_{n_j} }. \eqno (5.3)
$$
It follows from (4.1) that 
$$
\sum_{k = 1}^{m_j - n_j} \, \deg b_k  \le (M+2)(m_j - n_j) \eqno (5.4)
$$
and
$$
\deg s_{m_j}   = \sum_{k = 1}^{m_j} \, \deg b_k  
 = \deg q_{n_j}  + \sum_{k = 1}^{m_j - n_j} \, \deg b_k. \eqno (5.5)
$$
We infer from (5.3), (5.4) and (5.5) that
$$
\Vert s_{m_j}\Theta \Vert \le 
{1\over \vert s_{m_j}\vert \cdot2^{-2(M+2)(m_j - n_j)}} \cdot
\eqno (5.6)
$$

Since $\varphi$ is a non-increasing function
and $m_{j-1} + 1 = m_j - n_j$, we deduce from (5.1) that
$$
\eqalign{\varphi(\vert s_{m_j}\vert)\le 
\varphi(2^{{m_j}}) \le  2^{-2(M+2)(m_{j-1}+1)}
= 2^{-2(M+2)(m_j - n_j)}.\cr} \eqno (5.7)
$$
From  (5.2), (5.6) and (5.7), we thus obtain that
$$
|s_{m_j} | \cdot \Vert s_{m_j}  \, \Phi \Vert 
\cdot \Vert s_{m_j}  \, \Theta \Vert \le
\Vert s_{m_j}  \, \Theta \Vert \le
{1 \over |s_{m_j} | \cdot \varphi(|s_{m_j} |) } \cdot
$$
Since $j \ge 2$ is arbitrary, we have established that (2.1) holds.

It now remains to prove that $1$, $\Theta$ and $\Phi$ are independent 
over $\k[X]$. Therefore, we assume that they are dependent and 
we aim at deriving a contradiction. Let
$(A,B,C)$ be a non-zero triple of polynomials in $\k[X]$ satisfying 
$$
A \Theta + B \Phi + C=0. 
$$ 
Then, for any non-zero polynomial $q$ in $\k[X]$, we have
$$
\Vert q A\Theta \Vert=\Vert q B\Phi\Vert.
$$
In particular, we get
$$
\Vert s_{m_j} A\Theta\Vert= \Vert s_{m_j} B \Phi \Vert\leq \vert
B\vert\cdot 
\Vert s_{m_j}\Phi\Vert\ll |s_{m_j}|^{-1}, \eqno (5.8)
$$
for any $j \ge 2$. Here and below, the constants implied by $\ll$ depend
(at most) on $A$, $B$, $C$, $\Theta$ and $M$, but do
not depend on $j$.

On the other hand, we have constructed the sequence $(n_j)_{j\ge 1}$ 
in order to guarantee that  
$$
\vert s_{m_j}\Theta -s_{m_j-1}\vert\le  
{1\over |s_{m_j}| \cdot \varphi(|s_{m_j}|)}.\eqno (5.9)
$$
This implies that 
$$
\Vert s_{m_j}\Theta\Vert= \vert s_{m_j} \Theta - s_{m_j - 1}  \vert
$$
for $j$ jarge enough. 
Since by assumption the degree of
$b_{m_{j-1}+1}=t_{j-1}$ is either $M+1$ or $M+2$, 
we have $\deg b_{m_{j-1}+1}\not= \deg a_{n_j+1}$ 
and in particular $b_{m_{j-1}+1}\not=a_{n_j+1}$. 
Consequently, Lemma 3 implies that 
$$
\biggl\vert \Theta - {s_{m_j - 1} \over s_{m_j}} \biggr \vert
\ge   {1\over 2^{2(M+2)}\cdot  \vert q_{n_j}\vert^2}\gg 
{1\over 2^{2\deg q_{n_j}}},
$$
thus,
$$
\Vert s_{m_j} \Theta \Vert \gg 
2^{\deg s_{m_j}  - 2 \deg q_{n_j} }. \eqno (5.10)
$$
Moreover, we infer from (5.5) that
$\deg s_{m_j} \ge  \deg q_{n_j} + m_{j-1}$. Combined with (5.10),
this gives
$$
\vert s_{m_j}\vert\cdot \Vert s_{m_j}\Theta \Vert \gg
2^{2 m_{j-1}}. \eqno (5.11)
$$
For $j$ large enough, we deduce from (5.9) that 
$$
\vert s_{m_j} A \Theta- s_{m_j-1} A \vert < 2^{-1},
$$  
thus,
$$
\Vert s_{m_j}A\Theta\Vert=\vert s_{m_j}A\Theta- s_{m_j-1} A\vert\vert A\vert\cdot\Vert s_{m_j}\Theta\Vert.
$$
By (5.11), this yields
$$
\vert s_{m_j}\vert\cdot\Vert s_{m_j}A\Theta\Vert\gg 2^{2 m_{j-1}},
$$
which contradicts (5.8). 
This completes the proof of Theorem 1.  \cqfd

\vskip 8mm

\centerline{\bf 6. Proof of Theorems 2, 3 and 4}

\vskip 6mm

In all the proofs below we assume that $|\Theta| \le 1/2$ (replace 
$\Theta$ by $ 1/(X \Theta)$ if needed). 
The constants implied by $\ll$ may 
depend on $\Theta$ but not on $k$. 

\medskip

\noi {\bf Proof of Theorem 2.} 
Let $\Theta = [0,a_1, a_2, \ldots]$ and denote by $(p_n/q_n)_{n \ge 1}$
the sequence of its convergents.
The key observation for the proof of Theorem 2
is Lemma 1. Indeed, assume that the integer $n \ge 3$
is such that $a_1 \ldots a_n$ is a palindrome.
It then follows from Lemma 1 that
$q_n/q_{n-1}$ is very close to $1 / \Theta$. Precisely, we have
$$
\Vert q_{n-1} \Theta^{-1} \Vert \le 2^{-\deg q_{n-1}},
$$
by Lemma 2. Furthermore, (4.2) asserts that
$$
\Vert q_{n-1} \Theta \Vert = 2^{-\deg q_n}.
$$
Consequently, we get
$$
\vert q_{n-1}\vert^2\cdot  \Vert q_{n-1} \Theta \Vert
\cdot \Vert q_{n-1} \Theta^{-1} \Vert2^{-2 \deg q_{n-1}} \cdot \Vert q_{n-1} \Theta \Vert
\cdot \Vert q_{n-1} \Theta^{-1} \Vert < 1.
$$
This ends the proof. \cqfd

\goodbreak
\medskip

\noi {\bf Proof of Theorem 3.}
Assume now that $\Theta$ is in $\Bad$ 
and satisfies the assumption of Theorem 3.
Let $k \ge 1$ be an integer and let $P_k/Q_k$ be  
the last convergent to the rational number
$$
{P'_k \over Q'_k} := [0, U_k, V_k, \oU_k].
$$
Since, by assumption, $V_k$ is a palindrome, 
we obtain that the word $U_kV_k\oU_k$ is also a palindrome. 
Then, Lemma 1 implies that $P'_k = Q_k$.
Setting $r_k = |U_k|$ and $s_k=\vert V_k\vert$, 
we infer from Lemma 2 that
$$
\Vert Q_k \Theta \Vert\le  \, 2^{\deg Q_k} \,
2^{- 2 \deg q_{r_k +  s_k}}. \eqno (6.1)
$$
Observe that by Lemmas 1 and 2 we have
$$
\biggl| \Theta^{-1}- {Q'_k \over Q_k} \biggr|\ll 2^{- 2 \deg q_{r_k +  s_k}}
$$
and thus 
$$
\Vert Q_k\Theta^{-1}\Vert \ll 2^{\deg Q_k} \,
2^{- 2 \deg q_{r_k +  s_k}}. \eqno (6.2)
$$

Furthermore, it follows from (4.1) that 
$$ 
\deg Q_k =\deg q_{r_k} + \deg  q_{r_k+s_k}.
$$
Then, we get from (6.1) and (6.2) that 
$$
\vert Q_k\vert \cdot \Vert Q_k \Theta \Vert \cdot \Vert Q_k 
\Theta^{-1} \Vert\ll
2^{3 \deg q_{r_k}} \, 2^{- \deg q_{r_k+s_k}}. 
$$
In virtue of (2.3), this concludes the proof. 
 \cqfd

\medskip

\noi {\bf Proof of Theorem 4.}
Assume now that $\Theta$ and $\Phi$ satisfy the assumption of Theorem 4.
Let $k \ge 1$ be an integer and let $P_k/Q_k$ be  
the last convergent to the rational number
$$
{P'_k \over Q'_k} := [0, U_k, V^{[n_k]}].
$$
On the one hand, (4.2) gives  
$$
\Vert Q_k\Theta\Vert<{1\over \vert Q_k\vert}\cdot
$$
On the other hand, Lemma 1 implies that 
$$
{Q'_k \over Q_k}=[\oV^{[n_k]},\oU_k].
$$
Setting $r_k = |U_k|$ and $s_k=\vert V^{[n_k]}\vert$, 
we thus infer from Lemma 2 and (4.1) that
$$
\Vert Q_k \Phi \Vert \ll  \, 2^{Mr_k-ms_k} .
$$
In virtue of (2.4), this concludes the proof. 
 \cqfd
\goodbreak

\vskip 12mm

\centerline{\bf References}

\vskip 7mm

\vskip 4mm

\item{[1]}
B. Adamczewski and Y. Bugeaud,
{\it On the Littlewood conjecture in simultaneous
Diophantine approximation},
J. London Math. Soc. To appear.

\sm

\item{[2]}
J. V. Armitage,
{\it An analogue of a problem of Littlewood}, Mathematika 16 (1969),
101--105. Corrigendum and addendum, Mathematika 17 (1970),
173--178.

\sm

\item{[3]}
A. Baker,
{\it On an analogue of Littlewood's diophantine approximation
problem}, Michigan Math. J. 11 (1964), 247--250.

\sm

\item{[4]}
L. E. Baum and M.M. Sweet
{\it Continued fractions of algebraic power series in characteristic $2$},
Ann. of Math. 103 (1976), 593--610.

\sm

\item{[5]}
M. W. Buck and D. P. Robbins,
{\it The Continued Fraction Expansion of an Algebraic Power Series
Satisfying a Quartic Equation},
J. Number Theory 50 (1995), 335--344.

\sm

\item{[6]}
R. T. Bumby,
{\it On the analog of Littlewood's problem in power series fields},
Proc. Amer. Math. Soc. 18 (1967), 1125--1127.

\sm

\item{[7]} E. B. Burger, 
{\it On simultaneous Diophantine approximation 
in the vector space $\Q+ \Q\alpha$}, 
J. Number Theory 82 (2000), 12--24.

\sm

\item{[8]} 
J. W. S. Cassels and H. P. F. Swinnerton-Dyer, 
{\it On the product of three homogeneous linear forms 
and the indefinite ternary quadratic forms},  
Philos. Trans. Roy. Soc. London. Ser. A. 248 (1955), 73--96.

\sm

\item{[9]}
T. W. Cusick,
{\it Littlewood's Diophantine approximation problem for series}, $\!$
Proc. Amer. Math. Soc. 18 (1967), 920--924.

\sm

\item{[10]}
T. W. Cusick,
{\it Lower bound for a Diophantine approximation function},
Monasth. Math. 75 (1971), 398--401.

\sm

\item{[11]}
H. Davenport and D. J. Lewis, 
{\it An analogue of a problem of Littlewood}, Michigan Math. J. 
10 (1963) 157--160.

\sm

\item{[12]}
M. Einsiedler, A. Katok and E. Lindenstrauss,
{\it Invariant measures and the set of exceptions to the
Littlewood conjecture},
Ann. of Math. To appear.

\sm

\item{[13]}
T. Komatsu,
{\it Extension of Baker's analogue of Littlewood's Diophantine
approximation problem},
Kodai Math. J. 14 (1991), 335--340.

\sm

\item{[14]}
A. Lasjaunias,
{\it Diophantine Approximation and Continued Fraction Expansions 
of Algebraic Power Series in Positive Characteristic},
J. Number Theory 65 (1997), 206--225.

\sm

\item{[15]}
A. Lasjaunias,
{\it Quartic power series in $\F_3 ((T^{-1}))$ with bounded
partial quotients},
Acta Arith. 95 (2000), 49--59.

\sm
\item{[16]}
A. Lasjaunias,
{\it A survey of Diophantine approximation in fields of power series}, $\!$
Monatsh. Math. 130 (2000), 211--229.

\sm

\item{[17]}
J. E. Littlewood,
Some problems in real and complex analysis.
D. C. Heath and Co. Raytheon Education Co., 
Lexington, Mass., 1968.

\sm

\item{[18]}
K. Mahler,
{\it On a theorem of Liouville in fields of positive characteristic},
Canad. J. Math. 1 (1949), 397--400.

\sm

\item{[19]}
B. de Mathan,
{\it Quelques remarques sur la conjecture de Littlewood 
(en approximations diophantiennes simultan\'ees)},
S\'eminaire de Th\'eorie des Nombres, 1978--1979,
Exp. No. 5, 14 pp., CNRS, Talence, 1979.

\sm

\item{[20]}
B. de Mathan,
{\it Quelques remarques sur la conjecture de Littlewood 
(en approximations diophantiennes simultan\'ees)},
S\'eminaire de Th\'eorie des Nombres, 1980--1981,
Exp. No. 1, 20 pp., CNRS, Talence, 1981.

\sm

\item{[21]}
B. de Mathan,
{\it Simultaneous Diophantine approximation for algebraic functions
in positive characteristic},
Monatsh. Math. 111 (1991), 187--193.

\sm

\item{[22]}
B. de Mathan,
{\it A remark about Peck's method in positive characteristic},
Monatsh. Math. 133 (2001), 341--345.

\sm

\item{[23]}
W. H. Mills and D. P. Robbins,
{\it Continued fractions for certain algebraic power series},
J. Number Theory 23 (1986), 388--404.

\sm

\item{[24]} 
A. D. Pollington and S. Velani,
{\it On a problem in simultaneous Diophantine approximation:
Littlewood's conjecture},
Acta Math. 185 (2000), 287--306.

\sm

\item{[25]} 
M. Ratliff,
{\it The Thue--Siegel--Roth--Schmidt Theorem for algebraic functions},
J. Number Theory 10 (1978), 99--126.

\sm

\item{[26]} 
W. M. Schmidt,
{\it On continued fractions and Diophantine approximation in
power series fields},
Acta Arith. 95 (2000), 139--166.

\sm

\item{[27]}
D. S. Thakur,
Function field arithmetic,
World Scientific Publishing Co., Inc.,
River Edge, NJ, 2004.

\vskip15mm

\noindent Boris Adamczewski   \hfill{Yann Bugeaud}

\noindent   CNRS, Institut Girard Desargues 
\hfill{Universit\'e Louis Pasteur}

\noindent   Universit\'e Claude Bernard Lyon 1 
\hfill{U. F. R. de math\'ematiques}

\noindent   B\^at. Braconnier, 21 avenue Claude Bernard
 \hfill{7, rue Ren\'e Descartes}

\noindent   69622 VILLEURBANNE Cedex (FRANCE)   
\hfill{67084 STRASBOURG Cedex (FRANCE)}

\vskip2mm
 
\noindent {\tt Boris.Adamczewski@math.univ-lyon1.fr}
\hfill{{\tt bugeaud@math.u-strasbg.fr}}

\bye